\documentclass{article}
\usepackage{nips14submit_e,times}
\pdfoutput=1
\usepackage{standard}

\usepackage{times}
\usepackage{hyperref}
\usepackage{graphicx}
\usepackage{xypic}
\usepackage[all,2cell,ps]{xy}
\usepackage{tikz}
\usepackage{epsfig}
\usepackage{enumerate}
\usepackage{stmaryrd}
\usepackage{bbm}
\usepackage{algorithm}
\usepackage{algorithmic}
\usepackage{color}
\usepackage{paralist}

\newtheorem{defn}{Definition}

\newcounter{exmpl}

\newenvironment{mequation*}{$}{$}

\newcommand{\offdiag}{\text{off-diag}}
\newcommand{\Frob}{\text{Frob}}

\newcommand{\sqbb}[1]{\left\llbracket #1 \right\rrbracket}
\newcommand{\includepic}[2]{\includegraphics[width=#1\textwidth]{#2}}
\newcommand{\includepicc}[2]{\centerline{\includegraphics[width=#1\textwidth]{#2}}}

\renewcommand{\df}[1]{\textbf{#1}}

\makeatletter
\providecommand*{\cupdot}{\mathbin{\mathpalette\@cupdot{}}}
\newcommand*{\@cupdot}[2]{\ooalign{$\m@th#1\cup$\cr\hidewidth$\m@th#1\cdot$\hidewidth}}
\makeatother

\newcommand{\tikmxA}[1]{
\br{\,\begin{tikzpicture}[baseline=-13.5, scale=0.045]
\filldraw[gray] (0,0) rectangle +(#1,-#1); \end{tikzpicture}\,}}

\newcommand{\tikmxB}[2]{
\br{\,\begin{tikzpicture}[baseline=-13.5, scale=0.045] 
\filldraw[gray] (0,0) rectangle +(#1,-#1);
\foreach \i in {#1,...,#2}{
\filldraw[gray] (\i,-\i) rectangle +(1,-1);}
\end{tikzpicture}\,}}

\newcommand{\tikmxC}[2]{
\br{\,\begin{tikzpicture}[baseline=-13.5, scale=0.045]
\draw (0,0) rectangle +(17,-17);\draw (10,-9) node {#1}; \end{tikzpicture} \,}}

\title{Parallel MMF: a Multiresolution Approach to Matrix Computation}     
\author{Risi Kondor\m{{}^{\ast\dag}} Nedelina Teneva\m{{}^\ast} and Pramod K. Mudrakarta\m{{}^\ast}\\
\m{{}^\ast}Department of Computer Science, \m{{}^\dag}Department of Statistics\\
\texttt{\{risi,nteneva,pramodkm\}@cs.uchicago.edu}\\
The University of Chicago}

\nipsfinalcopy

\begin{document}
\maketitle
\vspace{-20pt}
\begin{abstract}
Multiresolution Matrix Factorization (MMF) was recently introduced as a method for 
finding multiscale structure and defining wavelets on graphs/matrices. 
In this paper we derive pMMF, a parallel algorithm for computing the MMF factorization. 
Empirically, the running time of pMMF scales linearly in the dimension for sparse matrices. 
We argue that this makes pMMF a valuable new computational primitive in its own right, 
and present experiments on using pMMF for two distinct purposes: 
compressing matrices and preconditioning large sparse linear systems. 
\end{abstract}

\section{Introduction}\label{sec: intro}

While the size of machine learning datasets continues to increase at the rate an order of magnitude or 
more every few years, the clock speed of commodity hardware has all but plateaued. 
Thus, increasingly, real world learning problems can only be tackled by algorithms that implicilty 
or explicitly exploit parallelism. 

Broadly speaking, there are two main approaches to parallelizing the 
optimization problems at the heart of machine learning algorithms. 
In the data parallel model, the dataset is divided into batches (shards) that are processed independently 
by separate processing units, but all processing units have (synchronous or asynchronous) access 
to a central resource that stores the current values of all the optimization variables, 
i.e., the model's parameters \cite{ho2013more}. 
In the parameter parallel approach, all processing units operate on the same data, 
but each one of them optimizes only a subset of the parameters \cite{parameter,lee2014model}. 
Data parallel algorithms often exploit the assumption that data subsets are \emph{i.i.d.} given the model parameters. 
However, such trivial parallelization is not suitable for model parallel algorithms, since it fails to capture the 
often non-trivial interactions between parameters and so taking independent parameter subsets could
lead to flawed estimates. On the other hand, model parallelization is often achieved by course graining the parameter space
(e.g., topological order or imposing some constraints on the graph induces by the parameter interaction), 
which might not be sufficient to discover the finer and more complex interactions between the parameters.

In contrast, in this paper we advocate a \emph{multiresolution} approach to parallelism, in which 
both the data and the parameters are parallelized, and communication 
between processing units is minimized. 
At the finest level of resolution, the data is divided into many shards, which are processed 
independently, determining the ``high frequency'' parameters of the model. These parameters are  
local to each shard, obviating the need for a central parameter server. 
Additionally, the data in each shard is compressed,  
so when it is redistributed across the cores for the second level, 
each shard becomes a compressed sketch of a larger subset of the original dataset, 
making it possible for the second level processing units to ``learn'' parameters that capture 
lower frequency, less localized features. Iterating this process leads to a multi-level 
algorithm that simultaneously uncovers the structure of the original data, and 
fits a model. 

The specific algorithm that we focus on in this paper is  
a parallelized version of the Multiresolution Matrix Factorization (MMF) process first described 
in \cite{MMFicml2014}. 
By itself, MMF is not a learning algorithm. 
However, parallel MMF is a gateway to efficiently performing a range of fundamental computational tasks, 
such as matrix compression, inversion, and solving linear systems. 
These tasks are critical building blocks of most learning algorithms. 
Due to space restrictions, the bulk of our numerical experiments, as well as the implementation 
details of our algorithm, which are critical to making pMMF scale to large problems, 
are relegated to the Supplement. 

\textbf{Notations.} In the following, \m{\sqb{n}} will denote the set \m{\cbrN{\oneton{n}}}. 
Given a matrix \m{A\tin\RR^{n\times n}} and two (ordered) sets \m{S_1,S_2\<\subseteq [n]},~\;  
\m{A_{S_1,S_2}} will denote the \m{\absN{\!S_1\!}\!\<\times\! \absN{\!S_2\!}} dimensional 
submatrix of \m{A} cut out by the rows indexed by \m{S_1} and the columns indexed by \m{S_2}. 
\m{\wbar{S_1}} will denote \m{[n]\!\setminus\! S_1}. 
\m{\ser{B}{\!\cupdot\!}{m}\<=[n]} denotes that the sets \m{\sseq{B}{m}} 
form a partition of \m{[n]}. \m{A_{:,i}} or \m{[A]_{:,i}} denotes the \m{i}'th column of \m{A}.  
\section{Parallel Multiresolution Matrix Factorization}\label{sec: MMF}

The Multiresolution Matrix Factorization (MMF) of a symmetric matrix 
\m{A\tin\RR^{n\times n}} is a multi-level factorization of the form 
\begin{equation}\label{eq: MMF1}
A\approx Q_1^\top\ldots Q_{L-1}^\top Q_L^\top H\, Q_L\ts Q_{L-1}\ldots Q_1,
\end{equation}
where \m{\sseq{Q}{L}} is a sequence of carefully chosen orthogonal matrices (rotations)  
obeying a number of constraints:
\begin{compactenum}[1.]
\item Each \m{Q_\ell} is chosen from some subclass \m{\Qcal} of highly sparse orthogonal matrices.  
In the simplest case, \m{\Qcal} is the class of \df{Givens rotations}, i.e., 
orthogonal matrices that only differ from the identity matrix in four matrix elements 
\begin{align*}
[Q_\ell]_{i,i}&=\cos \theta, &[Q_\ell]_{i,i}=-\sin \theta,\\
[Q_\ell]_{j,i}&=\sin \theta, &[Q_\ell]_{j,j}=\cos \theta,
\end{align*} 
for some pair of indices \m{(i,j)} and rotation angle \m{\theta}. 
Slightly more generally, \m{\Qcal} can 
be the class of so-called \m{\mathbf{k}}\df{--point rotations}, 
which rotate not just two, but \m{k} coordinates, \m{(\sseq{i}{k})}. 
\item The effective size of the rotations decreases according to a set schedule 
\m{n=\delta_0\geq \delta_1\geq\ldots\geq \delta_L}, i.e., there is a nested sequence 
of sets \m{[n]=S_0 \supseteq S_1\supseteq \ldots\supseteq S_L} with \m{\absN{\!S_\ell\!}=\delta_{\ell}} 
such that \m{[Q_\ell]_{\wbar{S_{\ell-1}},\wbar{S_{\ell-1}}}} 
is the \m{n\<-\delta_{\ell-1}} dimensional identity.  
\m{S_\ell} is called the \df{active set} at level \m{\ell}. 
In the simplest case, exactly one row/column is removed from the active set after each rotation. 
\item \m{H} is \m{\mathbf{S_L}}\df{--core-diagonal}, which means that it is all zero, except for (a) the submatrix 
\m{[H]_{S_L,S_L}} called the core and (b) the rest of the diagonal. 
\end{compactenum}
Moving the rotations in \rf{eq: MMF1} over onto the left hand side, the structure implied by 
the above conditions can be represented graphically as 
\begin{equation}\label{eq: MMF2}
\underset{\displaystyle Q_L^{\phantom{\top}}}{\tikmxB{5}{17}}
\hdots \underset{\displaystyle Q_2^{\phantom{\top}}}{\tikmxB{14}{17}}\nts
\underset{\displaystyle Q_1^{\phantom{\top}}}{\tikmxA{17}}
P\underset{\displaystyle A^{\phantom{\top}}}{\tikmxC{}{17}} P^\top \nts 
\underset{\displaystyle Q_1^\top}{\tikmxA{17}}\nts 
\underset{\displaystyle Q_2^\top}{\tikmxB{14}{17}}\hdots 
\underset{\displaystyle Q_L^\top}{\tikmxB{5}{17}}
\approx\underset{\displaystyle H^{\phantom{\top}}}{\tikmxB{3}{17}}\vspace{-15pt}
\end{equation}\vspace{5pt}\mbox{}\\ 
Here, for ease of visualization, \m{A} has been conjugated by a permutation matrix 
\m{P}, which ensures that \m{S_\ell\<=\cbrN{1,\ldots,\delta_\ell}} for each \m{\ell}.  
However, an actual MMF would not contain such an explicit permutation. 
In general, MMF is only an approximate factorization, because there is no guarantee that using a given 
number of rotations \m{A} can be brought into core-diagonal form with zero error. 
Most MMF factorization algorithms try to find \m{\sseq{Q}{L}} and \m{H} so as to minimize  
the squared Frobenius norm of the 
difference between the l.h.s and r.h.s of \rf{eq: MMF2}, called the \df{residual}. 

The original motivation for MMF in \cite{MMFicml2014} was to mimic 
the structure of fast orthogonal wavelet transforms. 
For example, when \m{A} is the Laplacian matrix of a graph \m{\Gcal} and 
\m{U=Q_L\ldots Q_2 Q_1}, the rows of \m{U} are interpreted as 
wavelets on the vertices of \m{\Gcal}. Those rows whose indices lie in 
\m{S_0\!\<\setminus\! S_1} are set by the first rotation matrix \m{Q_1}, 
and are not modified by the rest of the rotations. (If \m{Q_1} is a Givens rotation, 
there is only one such row, but for the more complicated rotations 
there might be several.)  
These rows are very sparse, and provide the lowest level, most local, highest frequency wavelets. 
The rows indexed by \m{S_1\!\<\setminus\! S_2} are determined by \m{Q_2 Q_1}, and correspond to 
level \m{2} wavelets, and so on. 
Writing the MMF as \m{A\approx U^\top\! H\tts U} suggests that the rows of \m{U} can 
also be interpreted as a hierarchically sparse PCA basis for \m{A}. 
Finally, since the size of the active set decreases after each rotation, 
defining \m{A_\ell\<=Q_\ell\ldots, Q_1 A\ts Q_1^\top \ldots Q_\ell^\top}, the sequence of transformations 
\begin{equation}\label{eq: transfseq}
A=A_0\mapsto A_1\mapsto A_2\mapsto \ldots \mapsto A_L\mapsto H
\end{equation}
is effectively a matrix compression scheme, which, according to \cite{MMFicml2014}, 
often significantly outperforms, for example, Nystr\"om methods. 
MMF is closely related to Diffusion Wavelets \cite{Coifman2004} and Treelets \cite{Lee2008}.

Due to the above properties, MMF is an attractive tool for uncovering the structure of large datasets. 
However, it has one fundamental limitation, which is its computational cost. 
The greedy factorization algorithm described in \cite{MMFicml2014} 
essentially follows the sequence of transformations in \rf{eq: transfseq}, at each level choosing 
\m{Q_\ell} and a set of rows/columns to be eliminated from the active set, 
so as to minimize their contribution to the final approximation error. 
Assuming the simplest case of each \m{Q_\ell} being a Givens rotation, this involves  
(a) finding the pair of vertices \m{(i,j)} involved in \m{Q_\ell}, and (b) finding the rotation 
angle \m{\theta}. In general, the latter is easy. However, 
finding the optimal choice of \m{(i,j)} (or, in the case of \m{k}--point rotations, 
\m{(\sseq{i}{k})}) is a combinatorial problem that scales poorly with \m{n}. 
The first obstacle is that the optimization is based on inner products between columns, so 
it requries computing the Gram matrix \m{G_\ell=A_{\ell-1}^\top\nts\nts A_{\ell-1}}  
at a complexity of \m{O(n^3)}. 
Note that this need only be done once: since rotations act on \m{G_\ell} the 
same way that they act on \m{A}, once we have \m{G_1\<=A^\top\! A}, each subsequent 
\m{G_\ell} can be efficiently derived via the recursion \m{G_{\ell+1}\!=Q_{\ell} G_{\ell} Q_{\ell}^\top}.  
The second obstacle is that searching for the optimal \m{(\sseq{i}{k})} has complexity \m{O(n^k)}.

The objective of the present paper is to construct a parallel MMF algorithm, \df{pMMF}, which, on 
typical matrices, assuming access to a sufficient number of processors, 
runs in time close to \emph{linear} in \m{n}. 
The ideas behind the new algorithm exploit the very structure in \m{A} that MMF factorizations 
pursue in the first place, namely locality at multiple levels of resolution. 

\subsection{Clustering}\label{sec: clustering}

The first and crucial step towards pMMF is to cluster the rows/columns of \m{A} into \m{m} clusters and 
only consider rotations between \m{k}-sets of rows/columns that belong to the same cluster. 
Letting \m{B_u} be the indices of the rows/columns belonging to cluster \m{u}, 
clustering has three immediate benefits: 
\begin{compactenum}[1.]
\item Instead of having to compute the full Gram matrix \m{G=A^\top\! A}, it is sufficient to 
compute the local Gram matrices \m{\cbrN{G^u\!\<=\!A_{B_u}^\top A_{B_u}}_{u=1}^m}. 
Assuming that the clustering is even, i.e., \m{\absN{\!B_u\!}\<=\Theta(c)} for some 
typical cluster size \m{c} (and therefore, \m{m\<=\Theta(n/c)}), 
this reduces the overall complexity of computing the Gram matrices from \m{O(n^3)} 
to \m{O(mc^2n)\<=O(c n^2)}. 
\item The complexity of the index search problem involved in finding each \m{Q_\ell} 
is reduced from \m{O(n^k)} to \m{O(c^k)}. In typical MMFs \m{\delta_L\nts\<=O(n)}, 
and the total number of rotations, \m{L},
scales linearly with \m{n}. Therefore, the total complexity of searching for rotations 
in the unclustered case is \m{O(n^{k+1})}, whereas with clustering it is \m{O(c^kn)}.    
\item The Gram matrices and the rotations of 
the different clusters are completely decoupled, therefore, 
on a machine with at least \m{m} cores, they can be computed in parallel, 
reducing the computation time of the above to \m{O(cn^2/m)\<=O(c^2n)} and 
\m{O(c^kn/m)\<=O(c^{k+1})}, respectively.  
\end{compactenum}
Using the clustering \m{B_1\<\cupdot B_2\<\cupdot\ldots\<\cupdot B_m\<=[n]} results in an MMF in which  
each of the \m{Q_\ell} matrices (and hence, also their product) are 
\m{(\sseq{B}{m})}--block-diagonal, as defined below.

\begin{defn}
Given a partition \m{B_1\<\cupdot B_2\<\cupdot\ldots\<\cupdot B_m} of \m{[n]}, we say that  
\m{M\tin\RR^{n\times n}} is \m{(\sseq{B}{m})}\df{--block-diagonal} if \m{M_{i,j}\<=0} unless  
\m{i} and \m{j} fall in the same cluster \m{B_u} for some \m{u}. 
\end{defn}


Clustering, in the form described above, decouples MMF into \m{m} completely independent subproblems. 
Besides the inherent instability of such an algorithm due to the vagaries of clustering algorithms, 
this approach is antithetical to the philosophy of multiresolution, since it cannot discover 
global features in the data: by definition, every wavelet will be local to only one cluster. 
The natural solution is to soften the clustering by repeatedly reclustering the data after a certain number of rotations. 
Writing out the full MMF again and grouping together rotations with the same clustering structure
\begin{equation*}
A\approx 
\underbrace{Q_1^\top\ldots Q_{l_1}^\top}_{\wbar{Q_1}^\top}
\underbrace{Q_{l_1\nts+1}^\top \ldots Q_{l_2}^\top}_{\wbar{Q_2}^\top} \ldots 
\underbrace{\ignore{Q_{l_{P-1}\nts+1}^\top} \ldots Q_{l_P}^\top}_{\wbar{Q_P}^\top} 
\;H\, 
\underbrace{Q_{l_{P}} \ldots \ignore{Q_{l_{P-1}+1}}}_{\wbar{Q_P}}\ldots  
\underbrace{Q_{l_2}\ldots Q_{l_1\nts+1}}_{\wbar{Q_2}}
\underbrace{Q_{l_1}\ldots Q_{1}}_{\wbar{Q_1}}
\end{equation*}
results in a factorization 
\begin{equation}
A\approx \wbar{Q}_1^\top \wbar{Q}_2^\top \ldots  \wbar{Q}_{P}^\top H\, \wbar{Q}_P\ldots  \wbar{Q}_{2} \wbar{Q}_1,
\end{equation}
where each \m{\wbar{Q}_p}, which we call a \df{stage}, is now a product of many \m{Q_\ell} elementary rotations, 
all conforming to the same block diagonal structure \m{(\sseq{B^p}{m_p})}. 
Note that the number of clusters, \m{m_p}, might not be the same across stages. 
In particular, since the active part of \m{A_\ell} progressively gets smaller and smaller,  
\m{m_p} will usually decrease with \m{p}. 
The top level pseudocode of pMMF, driven by this repeated clustering process, is given in 
Algorithm \ref{alg: top}. 

\begin{algorithm}[t]
\begin{algorithmic}
\STATE \textbf{Input:} a symmetric matrix \m{A\tin\RR^{n\times n}}
\STATE \textbf{Set} \m{A_0\leftarrow A}
\STATE \textbf{for} (\m{p\<=1} to \m{P})\,\m{\{}
\STATE ~~\textbf{cluster} \,the active columns of \m{A_{p-1}} to get \m{(B^p_1,\ldots,B^p_m)}
\STATE ~~\textbf{reblock} \m{\, A_{p-1}} according to \m{(B^p_1,\ldots,B^p_m)}
\STATE ~~\textbf{for} (\m{u\<=1} to \m{m})~~\m{\wbar{Q}_{p,u}\leftarrow} 
\textbf{FindRotationsInCluster(}\m{p,u}\textbf{)}
\STATE ~~\textbf{for}\,\m(\m{u\<=1}~\textbf{to}~\m{m})~\,\textbf{for}\,(\m{v\<=1}~\textbf{to}~\m{m})~\;\textbf{set}~    
\m{\sqbb{A_{p}}_{u,v}\nts\leftarrow {\wbar{Q}_{p,u}}\sqbb{A_{p-1}}_{u,v} \wbar{Q}_{p,v}{\!\!\!\!\!}^{\top}} 
\STATE ~~\textbf{merge} \m{(\wbar Q_{p,1},\ldots,\wbar Q_{p,m})} \textbf{into} \m{\wbar{Q}_p}
\STATE \}
\STATE \m{H\leftarrow} the core of \m{A_L} plus its diagonal 
\STATE \textbf{Output:} \m{(H,\wbar{Q}_1,\ldots,\wbar{Q}_p)}
\end{algorithmic}
\caption{\label{alg: top}~\textbf{pMMF}~\;(top level of the algorithm)}
\end{algorithm}

\subsection{Randomized greedy search for rotations}

The second computational bottleneck in MMF 
is finding the \m{k} rows/columns involved in each rotation. 
To address this, we use a randomized strategy, whereby first a single row/column \m{i_1} 
is chosen from the active set (within a given cluster) unformly at random, and then \m{k\<-1} further  
rows/columns \m{i_2,\ldots,i_k} are selected from the same cluster 
according to some separable objective function \m{\phi(i_2,\ldots,i_{k})} related to minimizing the 
contribution to the final error. For simplicity, in pMMF we use 
\[\phi(i_2,\ldots,i_{k})=\sum_{r=2}^{k} 
\fr{\inp{\tts[A_{\ell-1}]_{:,i_1},[A_{\ell-1}]_{:,i_r}}}{\nmN{[A_{\ell-1}]_{:,i_r}}},\]
i.e., \m{[A_{\ell-1}]_{:,i_1}} is rotated with the \m{k\<-1} other columns that it has 
the highest normalized inner product with in absolute value. 
Similarly to \cite{MMFicml2014}, the actual rotation angle 
(or, in the case of \m{k}'th order rotations, the \m{k\<\times k} non-trivial submatrix of \m{Q_\ell}) 
is determined by diagonalizing \m{[G_\ell]_{(i_1\ldots i_k),(i_1\ldots i_k)}} at a cost of only \m{O(k^3)}. 
This aggressive randomized-greedy strategy reduces the complexity of finding each rotation to \m{O(c)}, 
and in our experience does almost as well as exhaustive search. 
The criterion for elimination is minimal off-diagonal norm, 
\m{\nmN{\!A_{:,i}\!}_{\offdiag}=(\nmN{\!A_{:,i}\!}^2-A_{i,i}^2)^{1/2}}, 
because \m{2\ts\nmN{\!A_{:,i}\!}_{\offdiag}^2} is the contribution of eliminating row/column \m{i} to the final error. 

\begin{algorithm}[b]
\begin{algorithmic}
\STATE \textbf{Input:} a matrix \m{\mathcal{A}\<=[A_p]_{:,B_u}\tin \RR^{n\times c}} made up of 
the \m{c} columns of \m{A_{p-1}} forming cluster \m{u} in \m{A_p}\\ 
\STATE \textbf{compute} the Gram matrix \m{G\<=\Acal^\top\! \Acal}\\
\STATE \textbf{set}\; \m{I\<=[c]}~~ (the active set) \\
\STATE \textbf{for}~(\m{s\<=1}~\textbf{to}~\m{\lfloor \eta c\rfloor})\m{\{}\\
\STATE ~~\textbf{select}~ \m{i\tin I} uniformly at random\\ 
\STATE ~~\textbf{find}~ \m{j\<=\argmax_{I\setminus\cbr{i}} 
\abs{\inpN{\Acal_{:,i},\Acal_{:,j}}}/\nmN{\nts \Acal_{:,j}\nts}}\\ 
\STATE ~~\textbf{find}~\;the Givens rotation~\m{q_s}~of columns \m{(i,j)} as described in the text\\ 
\STATE ~~\textbf{set}~\m{\Acal\leftarrow q_s \Acal\ts q_s^\top}\\
\STATE ~~\textbf{set}~\m{G\leftarrow q_s\tts G\tts q_s^\top}\\
\STATE ~~\textbf{if}~\:\m{\nm{\!\Acal_{i,:}}_{\offdiag}\!\<< \nm{\!\Acal_{j,:}}_{\offdiag}}~\textbf{eliminate}~\m{i}
\textbf{;~ otherwise eliminate}~\m{j}\\
\STATE \m{\}}
\STATE \textbf{Output:} \m{\wbar Q_{p,u}=q_{\lfloor \eta c\rfloor}\ldots q_2\tts q_1}
\end{algorithmic}
\caption{~\;\textbf{FindRotationsInCluster(}\m{p,u}\textbf{)} --- here \m{k\<=2} and 
\m{\eta} is the compression ratio
\label{alg: greedy}}
\end{algorithm} 

\subsection{Blocked matrices}

In a given cluster \m{u} of a given stage \m{p}, the local Gram matrix \m{G^u}, and 
the rotations can be determined from the columns belonging to just that cluster. 
However, subsequently, these rotations need to be applied to the entire matrix, from both the right 
and the left, which cuts across clusters. 
To be able to perform this part of the algorithm in parallel as well, we partition 
\m{A} not just column-wise, but also row-wise. The resulting data structure is called a symmetrically 
blocked matrix (c.f., \cite{buluc2012}). 

\begin{defn}
Given a matrix \m{A\tin\RR^{n\times n}} and a partition \m{B_1\<\cupdot B_2\<\cupdot\ldots\<\cupdot B_m} of \m{n}, 
the \m{(u,v)} block of \m{A} is the submatrix \m{\sqbb{A}_{u,v}\!:=A_{B_u,B_v}}.   
The \df{symmetric blocked matrix} form of \m{A} consists of the \m{m^2} 
separate matrices \m{\cbrN{\sqbb{A}_{u,v}}_{u,v=1}^{m}}. 
\end{defn}

In pMMF, the matrix \m{A_\ell} is always maintained in blocked matrix form, where the block structure 
is dictated by the clustering of the current stage. 
For large matrices, the individual blocks can be stored on separate cores or separate machines, and 
all operations, including computing the Gram matrices, are performed in a block-parallel fashion. 
This further reduces the time complexity of the 
Gram matrix computation from \m{O(nc^2)} to \m{O(c^3)}. 
Assuming \m{m_p^2}--fold parallelism, and a total of \m{\eta c} rotations in stage \m{p}, 
the overall time needed to apply all of these rotations to the entire matrix scales with 
\m{O(\eta kc^2)}. 

\begin{figure}
\centering
\begin{minipage}{.4\textwidth}
\tiny 
\setlength{\tabcolsep}{5pt}
\begin{tabular}{|p{9pt}|p{9pt}|p{9pt}|p{9pt}|p{9pt}|}
\hline
\m{\!\!\!\sqbb{\nts M\nts }_{11}^{\phantom{M^M}}}&\m{\!\!\!\sqbb{\nts M\nts}_{12}}&\m{\!\!\!\sqbb{\nts M\nts}_{13}}&\m{\!\!\!\sqbb{\nts M\nts}_{14}}&\m{\!\!\!\sqbb{\nts M\nts}_{15}}\\[4pt]
\hline
\m{\!\!\!\sqbb{\nts M\nts }_{21}^{\phantom{M^M}}}&\m{\!\!\!\sqbb{\nts M\nts}_{22}}&\m{\!\!\!\sqbb{\nts M\nts}_{23}}&\m{\!\!\!\sqbb{\nts M\nts}_{24}}&\m{\!\!\!\sqbb{\nts M\nts}_{25}}\\[4pt]
\hline
\m{\!\!\!\sqbb{\nts M\nts }_{31}^{\phantom{M^M}}}&\m{\!\!\!\sqbb{\nts M\nts}_{32}}&\m{\!\!\!\sqbb{\nts M\nts}_{33}}&\m{\!\!\!\sqbb{\nts M\nts}_{34}}&\m{\!\!\!\sqbb{\nts M\nts}_{35}}\\[4pt]
\hline
\m{\!\!\!\sqbb{\nts M\nts }_{41}^{\phantom{M^M}}}&\m{\!\!\!\sqbb{\nts M\nts}_{42}}&\m{\!\!\!\sqbb{\nts M\nts}_{43}}&\m{\!\!\!\sqbb{\nts M\nts}_{44}}&\m{\!\!\!\sqbb{\nts M\nts}_{45}}\\[4pt]
\hline
\m{\!\!\!\sqbb{\nts M\nts }_{51}^{\phantom{M^M}}}&\m{\!\!\!\sqbb{\nts M\nts}_{52}}&\m{\!\!\!\sqbb{\nts M\nts}_{53}}&\m{\!\!\!\sqbb{\nts M\nts}_{54}}&\m{\!\!\!\sqbb{\nts M\nts}_{55}}\\[4pt]
\hline
\end{tabular}
\end{minipage}
\begin{minipage}{.29\textwidth}
\centerline{\includegraphics[width=.65\textwidth, angle=90]{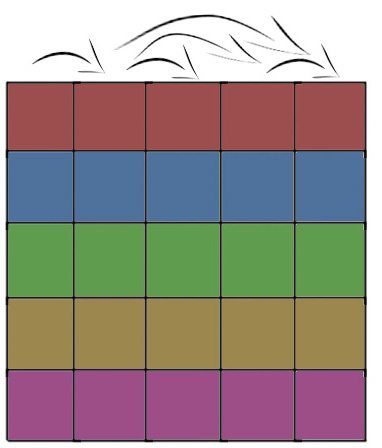}}
\end{minipage}
\begin{minipage}{.29\textwidth}
\includepicc{.65}{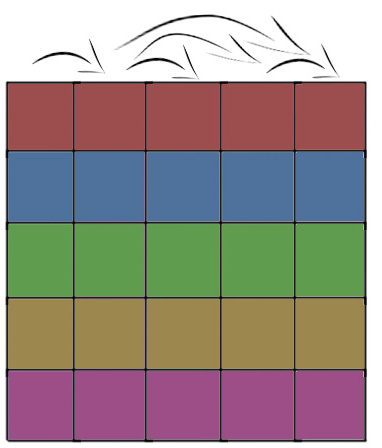}
\end{minipage}
\vspace{-6pt}
\caption{\label{fig: blocked} Schematic of a blocked matrix \m{M} with \m{5\<\times 5} blocks. 
For the sake of visual clarity, we assumed that the blocks are contiguous, 
but, in general, this is not the case.  
The reblocking process involves first reorganizing the rows acoording to the new structure, then 
reorganizing the columns. To perform this efficiently, the first 
operation is done in parallel for each column of blocks of the original matrix, and 
the second operation is done in parallel for each row of blocks. 
}
\end{figure}

The blocked matrix data structure is ideally suited to carrying out each stage of MMF 
on a parallel system, 
because (except for summary statistics) no data needs to be communicated 
between the different blocks. 
However, changing the block structure of the matrix from one clustering to another 
can incur a large communication overhead. 
To retain \m{m}--fold parallelism, the reblocking is carried out in two phases: 
first, each column of blocks is reblocked row-wise, then each row of blocks 
in the resulting new blocked matrix is reblocked column-wise (Figure \ref{fig: blocked}). 

\subsection{Sparsity and matrix-free MMF arithmetic}\label{sec: sparsity}

Ultimately, pMMF is intended for factoring matrices that are sparse, but whose dimensionality 
is in the hundreds of thousands or millions. As the factorization progresses, the fill-in 
(fraction of non-zeros) in \m{A_p} will increase, but at the same time,  
the active part of \m{A_\ell} will progressively shrink. 
This means that in practice, given a sufficiently highly parallel system, 
the overall complexity can still scale roughly linearly with the number of non-zeros in \m{A}. 

The complete factorization appearing on the r.h.s.\;of \rf{eq: MMF1} we denote \m{\tilde A}. 
Storing \m{\tilde A} by storing \m{H} and the \m{\cbr{Q_\ell}} matrices  
separately, the space complexity scales roughly linearly in the dimension. 
On the other hand, computing \m{\tilde A} explicitly as a dense \m{\RR^{n\times n}} 
matrix is usually unfeasible. 
Therefore, when applying the computed factorization, for example, as a preconditioner 
(Section \ref{sec: preconditioning}),   
which requires repeatedly multilying a vector \m{v} by \m{\tilde A}, we use the so-called matrix-free approach: 
\m{v} is stored in the same blocked form as \m{A}, the rotations are applied 
individually, and as the different stages are applied to \m{v}, the vector goes through an analogous 
reblocking process to that described for \m{A}. 
The complexity of matrix-free MMF/vector multiplication is \m{O(kpn)}. 
Inverting \m{\tilde A}, which is also critical for downstream applications, involves 
inverting the entries on the diagonal of \m{H} and inverting the core matrix, thus 
the overall complexity of MMF inversion is \m{O(n+\delta_L^3)}.

The theoretical complexity of the main components of pMMF are summarized in Table \ref{tbl: complexity}. 
Of course, requiring \m{m^2}--fold parallelism as \m{m\<\to\infty} is an abstraction. 
Note, however, that even the total operation count scales with \m{\gamma n^2}, which is just the 
number of non-zeros in the original matrix. 
The plots in Figure \ref{fig: time1} and similar plots in the Supplement confirm that on many real 
world datasets, particularly, matrices coming from sparse network graphs, the wall clock time 
of pMMF tends to scale linearly with the dimension. Also note that in these experiments \m{n} is on the 
order of \m{10^4\sim 10^5}, yet the factorization time is on the order of just one minute. 

\begin{figure}[t]
\centering
\includepic{.24}{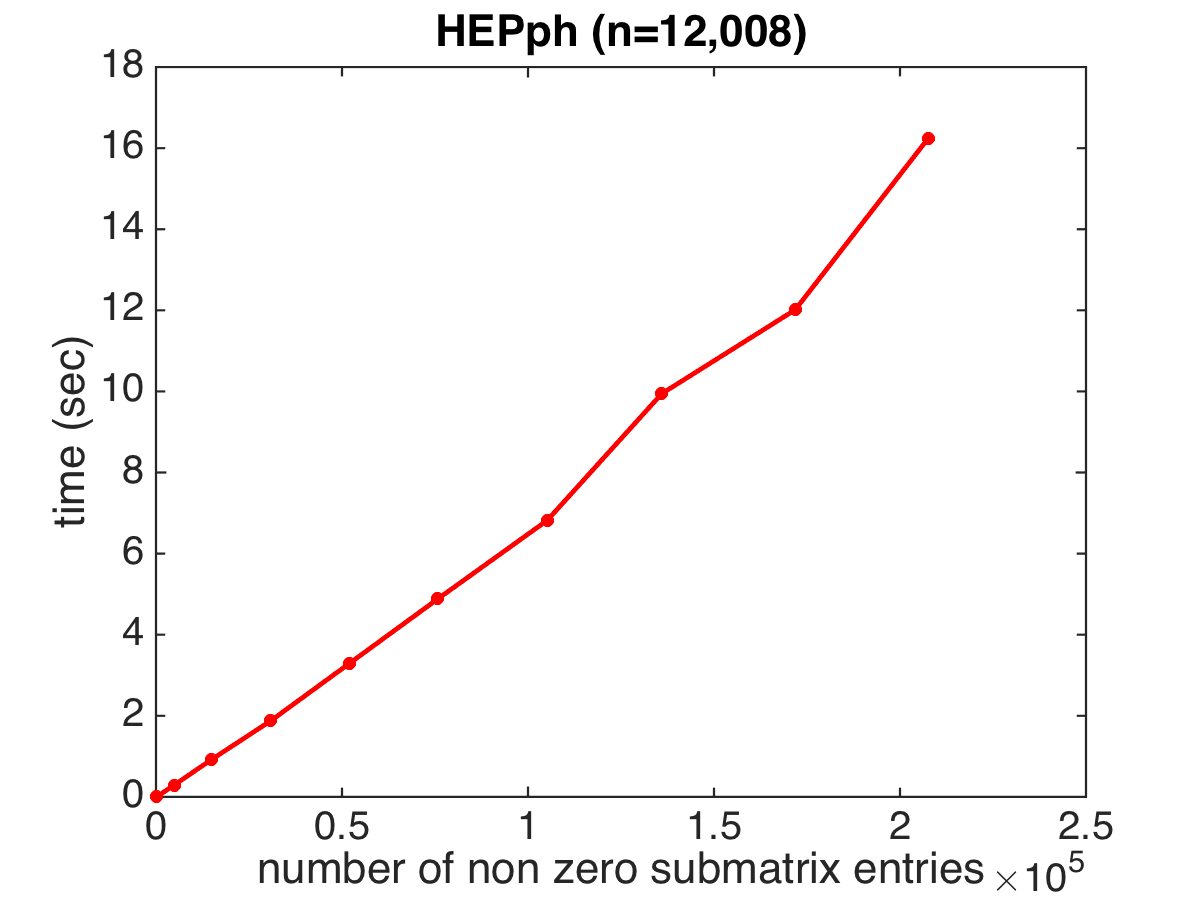}
\includepic{.24}{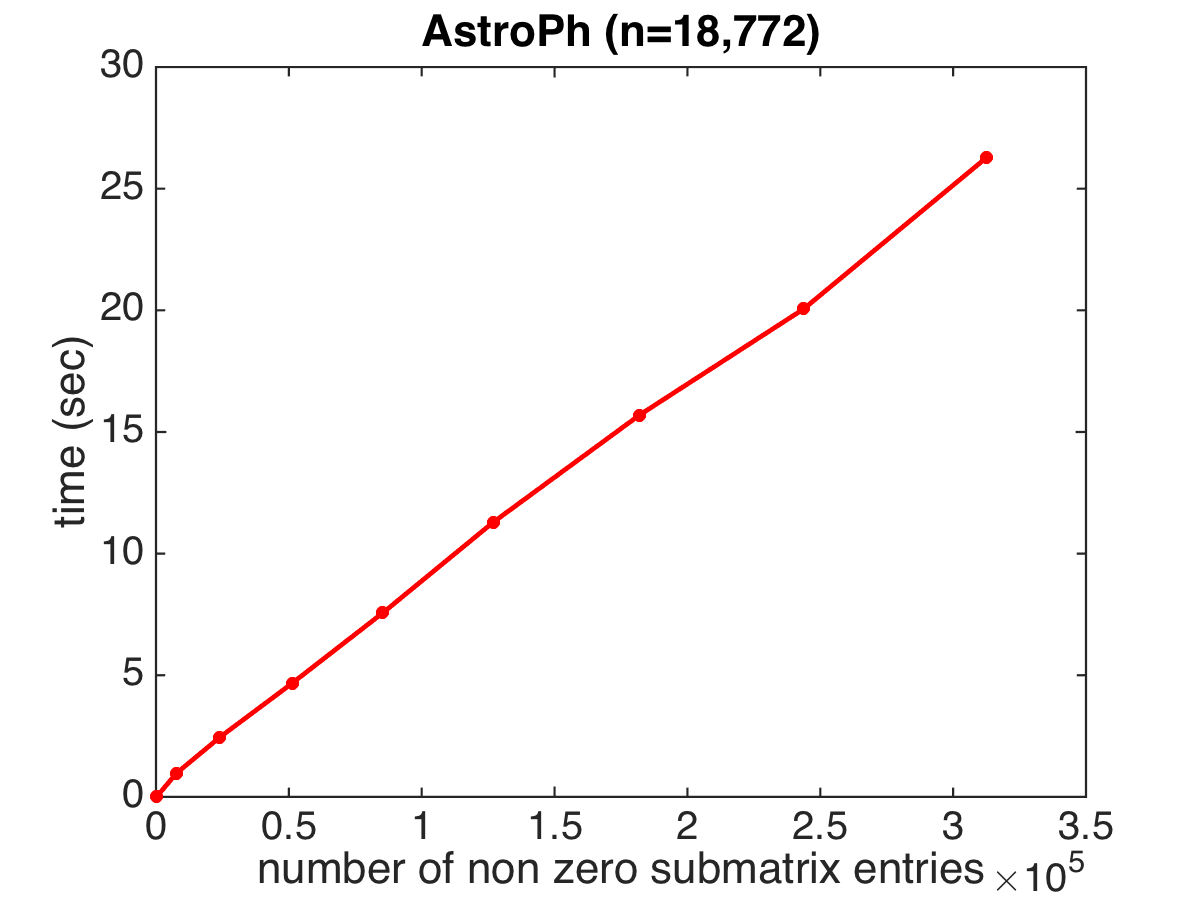}
\includepic{.24}{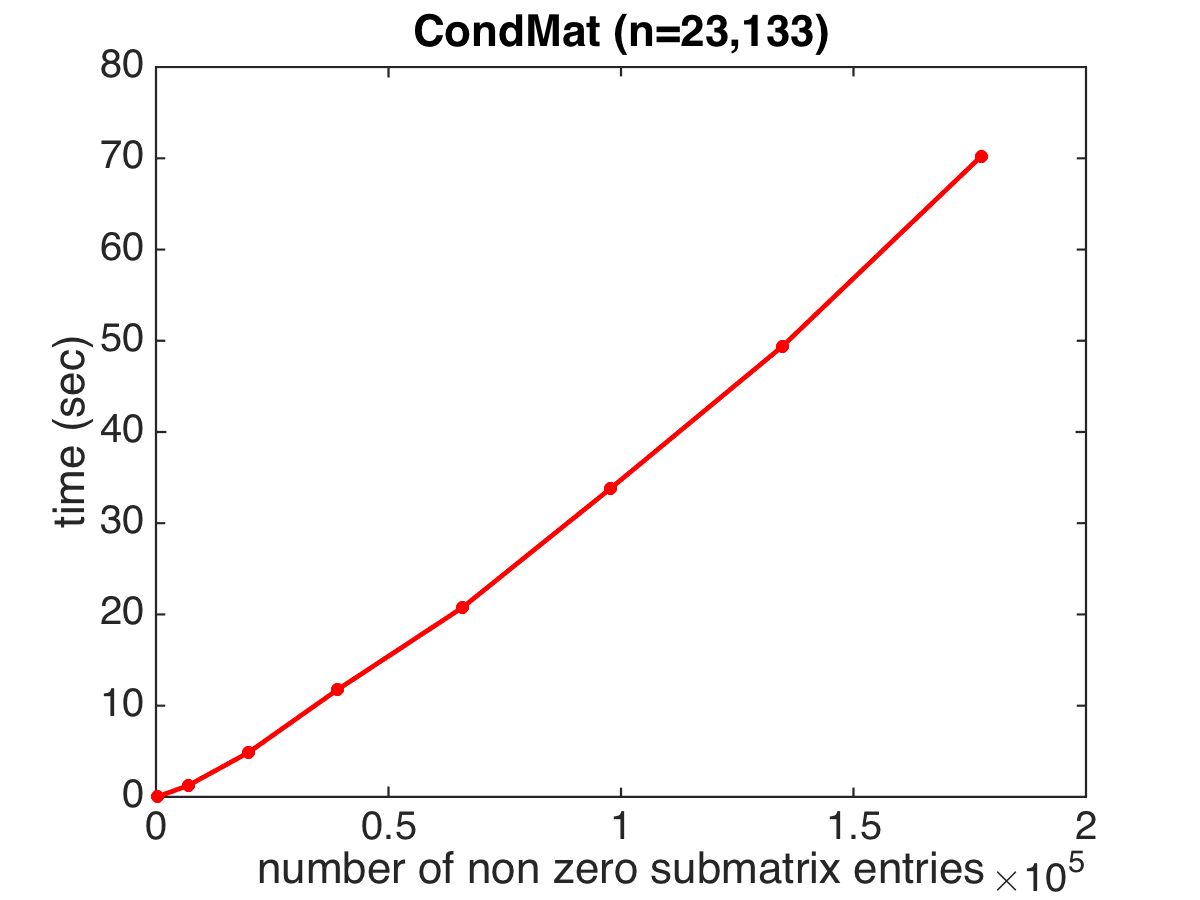}
\includepic{.24}{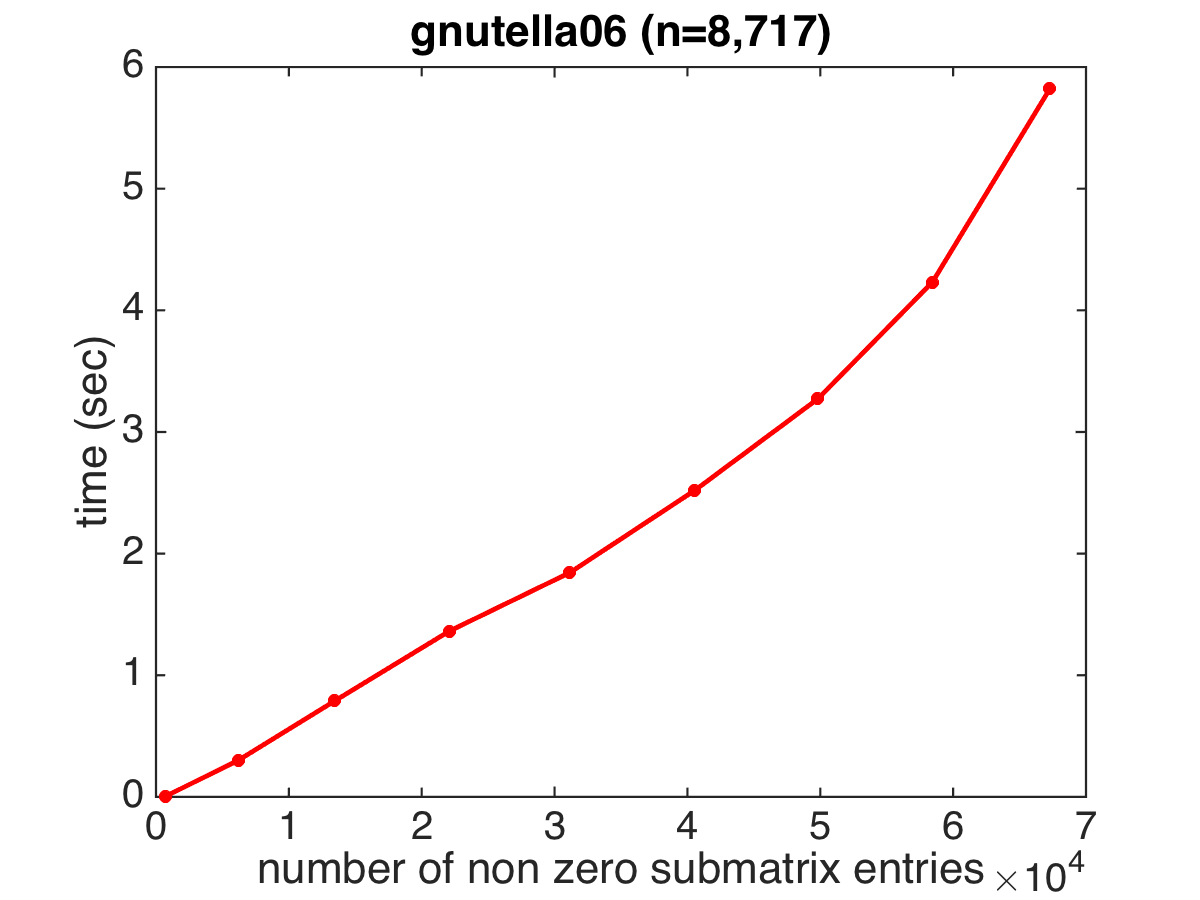}
\vspace{-7pt}
\caption{\label{fig: time1}  
Wall clock time of pMMF in seconds as a function of the number of non-zeros, \m{\gamma n^2}, in \m{A} on a 32 core 2.6 GHz machine. 
\m{A} is derived from the Laplacian of standard benchmark sparse graphs (see Supplement), taking 
submatrices of different sizes.
The figures confirm that in practice pMMF often scales linearly in the 
number of nonzeros, hence, for bounded degree graphs, also in \m{n}.  
}
\end{figure}

\ignore{
Given a set \m{\Xcal} and a space \m{L(\Xcal)} of sufficiently regular functions 
on \m{\Xcal}, \df{multiresolution analysis} \cite{Mallat1989} is a way of stratifying \m{L(\Xcal)} 
by constructing a nested sequence of subspaces  
\begin{equation}\label{eq: MRA1}
L(\Xcal)\supset \ldots\supset V_{-1}\supset V_0\supset V_1\supset V_2\ldots
\end{equation}
comprising functions at increasing levels of smoothness but decreasing level of detail. 
In the most classical setting \m{\Xcal} is \m{\RR^d}, and smoothness 
is defined in terms of a quantity such as 
\m{\eta=(\inp{f,\Delta f}/\inp{f,f})^{-1}}, where \m{\Delta=\dd^2/\dd x_1^2+\ldots+\dd^2/\dd x_d^2}, 
is the Laplace operator, and \m{\inp{f,g}=\int\! f(x)\tts g(x)\ts dx}. 
The way to look at \rf{eq: MRA1} is that 
each stage \m{V_\ell\to V_{\ell+1}} splits \m{V_\ell} into an orthogonal 
sum \m{V_{\ell+1}\<\oplus W_{\ell+1}}, where \m{W_\ell} is the space that captures detail at level \m{\ell+1}, 
while \m{V_{\ell+1}} retains the part of signals at all lower levels of resolution. 
To turn this recursive splitting process into a computational tool, one needs to 
define an orthogonal basis \m{\cbrN{\phi_{\ell,m}}{}_{m=1}^{d_\ell}} for each \m{V_\ell}, 
and a similar orthogonal basis \m{\cbrN{\psi_{\ell,m}}{}_{m=1}^{\delta_{\ell-1}-\delta_\ell}} for each \m{W_\ell}. 
It is the \m{\psi_{\ell,m}} basis functions that are reffered to as \df{wavelets}. 
Specifically, by the \df{wavelet transform} of a function \m{f\tin V_0} 
up to level \m{L}, we mean an expansion of the form 
\[f=\sum_{m=1}^{\delta_L}\alpha_m\,\phi_{L,m}+
\sum_{\ell=1}^{L}\sum_{m=1}^{\delta_\ell-\delta_{\ell-1}} \beta_{\ell,m}\,\psi_{\ell,m},\]
which expresses \m{f} in terms of the union of the bases of \m{\seq{W}{L}} and \m{V_L}. 
The key to making such wavelet transforms fast is to ensure that each constitutent  
\m{\cbrN{\phi_{\ell,m}}_m\mapsto \cbrN{\phi_{\ell}}\cup \cbrN{\psi_{\ell,m}}_m} transform 
is very sparse.  

In a recent paper, Kondor et al. \cite{Kondor} defined the Multiresolution 

The resulting theory leads to wavelets, which have proved themselves to 
be one of the most useful and versatile tools in all of applied mathematics. 
In recent years, the focus has been on extending multiresolution analysis, 
together with the corresponding 
to less structured settings, such as when \m{\Xcal} is a graph, 
and exploiting the properties 
of the resulting wavelet transforms. 

Regardless  
}
\section{pMMF Compression}\label{sec: compression}

Most, if not all, machine learning algorithms reduce to linear algebra operations or optimization over 
large instance/feature matrices or instance/instance similarity matrices. 
The classical example is, of course, kernel methods, which reduce to convex optimization 
involving the so-called Gram matrix, a symmetric positive semi-definite (p.s.d.) matrix of size \m{n\<\times n}, where 
\m{n} is the number of training examples. 
Despite their many attractive properties and their large literature, the applicability 
of kernel methods to today's large datasets is limited by the fact that ``out of the box'' 
their computation time scales with about \m{n^3}. 
This issue (not just for kernel methods, but more broadly) 
has catalyzed an entire area focused on compressing or ``sketching'' matrices 
with minimal loss. 

In the symmetric (and p.s.d) case, most sketching algorithms 
approximate \m{A\tin\RR^{n\times n}} in the form \m{\tilde A\<=C W^\dag C^\top}, where \m{C} is a judiciously chosen 
\m{\RR^{n\times m}} matrix with \m{m\<\ll n}, and \m{W^\dag} is computed by taking the pseudo-inverse 
of a certain matrix that is of size only \m{m\times m}.  
The algorithms mainly differ in how they define \m{C}:
\begin{compactenum}[(i)]
\item Projection based methods set each column of \m{C} to be a random linear combination of the columns 
of the original matrix \m{A}, and use Johnson--Lindenstrauss type arguments to show that the 
resulting low dimensional random sketch preserves most of the information at least about the part of \m{A} 
spanned by its high eigenvalue eigenvectors \cite{Halko2011}. 
These methods come with strong guarantees, but suffer from the cost of  
having to compute \m{m} dense linear combinations of the columns of \m{A}. 
Even if \m{A} was sparse, this process destroys the sparsity.  
\item Structured projections are a twist on the above idea, replacing the random projection 
with a fixed, dense, but efficiently computable, basis transformation (and subsampling in that basis), 
such as the fast Hadamard transform or the fast Fourier transform \cite{Ailon2009,fastfood2013}. 
\item In contrast to (i) and (ii), Nystr\"om algorithms construct \m{C} by choosing a certain number of 
\emph{actual} columns of \m{A}, usually by random sampling.  
Here, the focus has shifted from uniform sampling \cite{Williams2001,Fowlkes2004}, 
via \m{l_2}--norm sampling \cite{Drineas2006all}, to sampling based on so-called leverage scores  
\cite{Mahoney2009CUR, MahoneyRandomized, Gittens2013}, which, for low rank matrices, can be shown 
to be optimal. Further recent developments include the ensemble Nystr\"om method \cite{Kumar2009} 
and the clustered Nystr\"om algorithm \cite{Zhang2010}. 
Finally, a number of adaptive algorithms have also been proposed \cite{Deshpande2006,Kumar2012,Wang2013}. 
\end{compactenum}

pMMF can also be regarded as a sketching method, in the sense that \rf{eq: MMF2} compresses \m{A} into 
an \m{m\!:=\!\delta_L\nts\<\ll n} dimensional core via a very fast series of orthogonal 
transforms \m{Q_L\ldots Q_2 Q_1}, which play a role analogous to \m{C^\top}\!\!. 
In contrast to the other methods, however, MMF also retains the entries on the diagonal of \m{H}.   
From the point of view of downstream computation, this incurs very little extra cost. 
For example, if the purpose of sketching \m{A} is to compute its inverse, 
maybe for use in a Gaussian process or ridge regression, 
then \m{\tilde A^{-1}} can be easily computed by inverting the core of \m{H} (as in the Nystr\"om methods), 
and just taking the inverse of all other matrix elements on the diagonal. 

\begin{figure}[t]
\centering
\begin{minipage}{.245\textwidth}
\centerline{\includegraphics[width=1.12\textwidth]{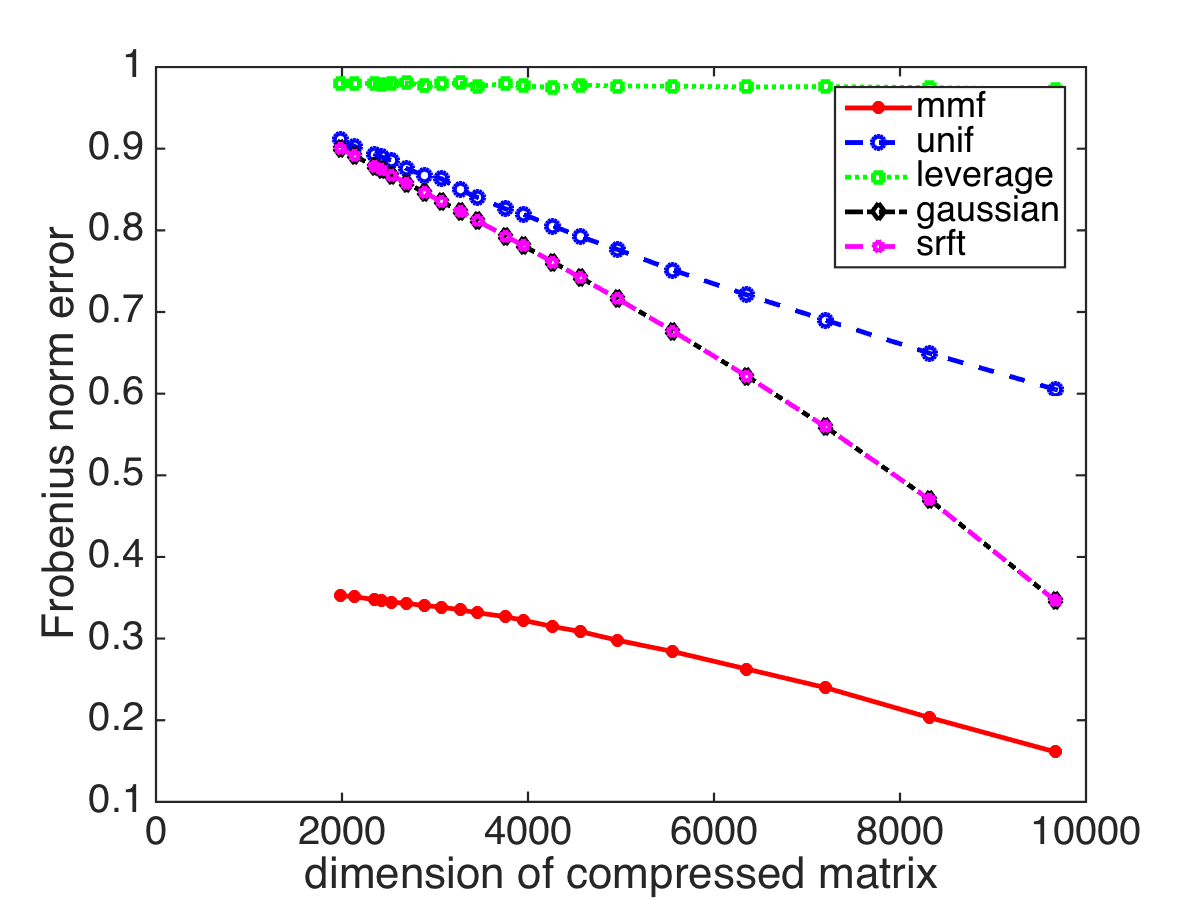}}
\end{minipage}
\begin{minipage}{.245\textwidth}
\centerline{\includegraphics[width=1.12\textwidth]{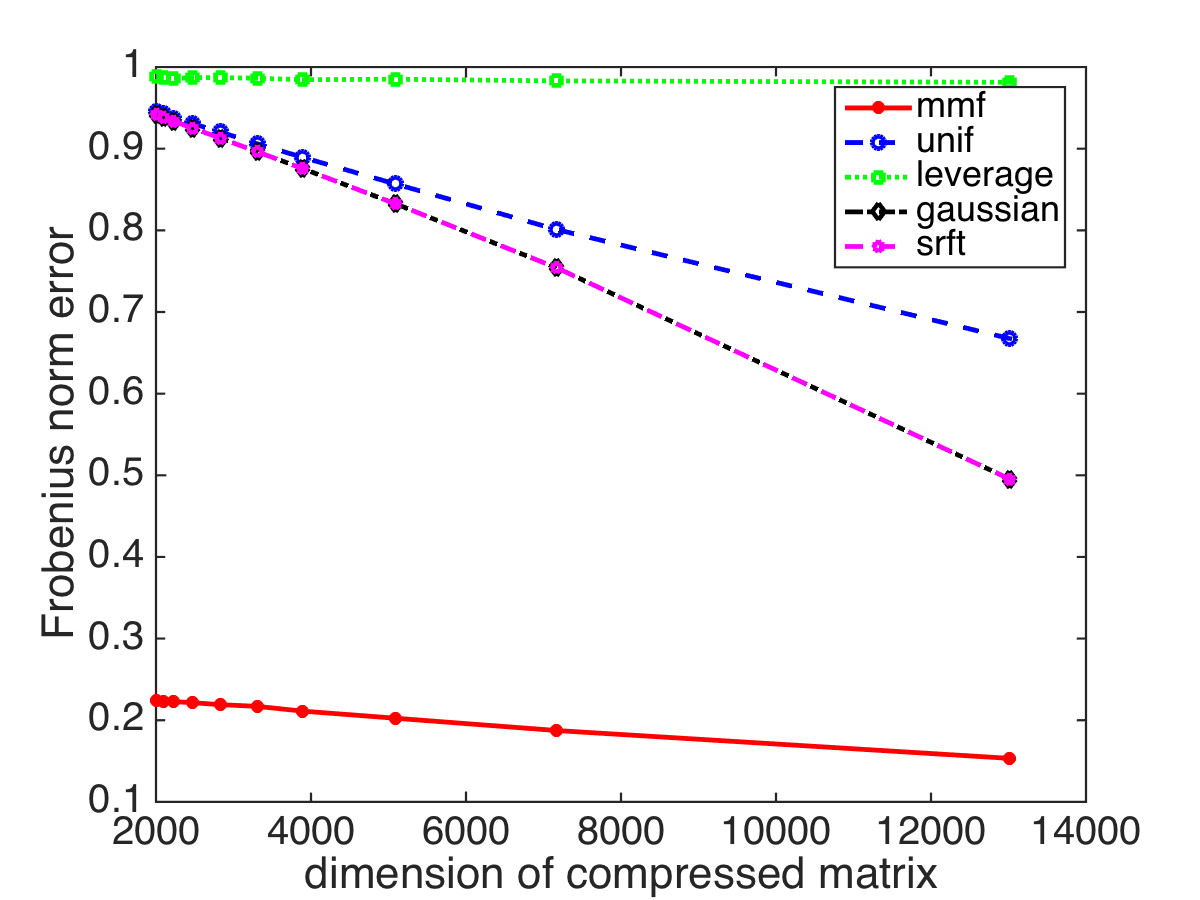}}
\end{minipage}
\begin{minipage}{.245\textwidth}
\centerline{\includegraphics[width=1.12\textwidth]{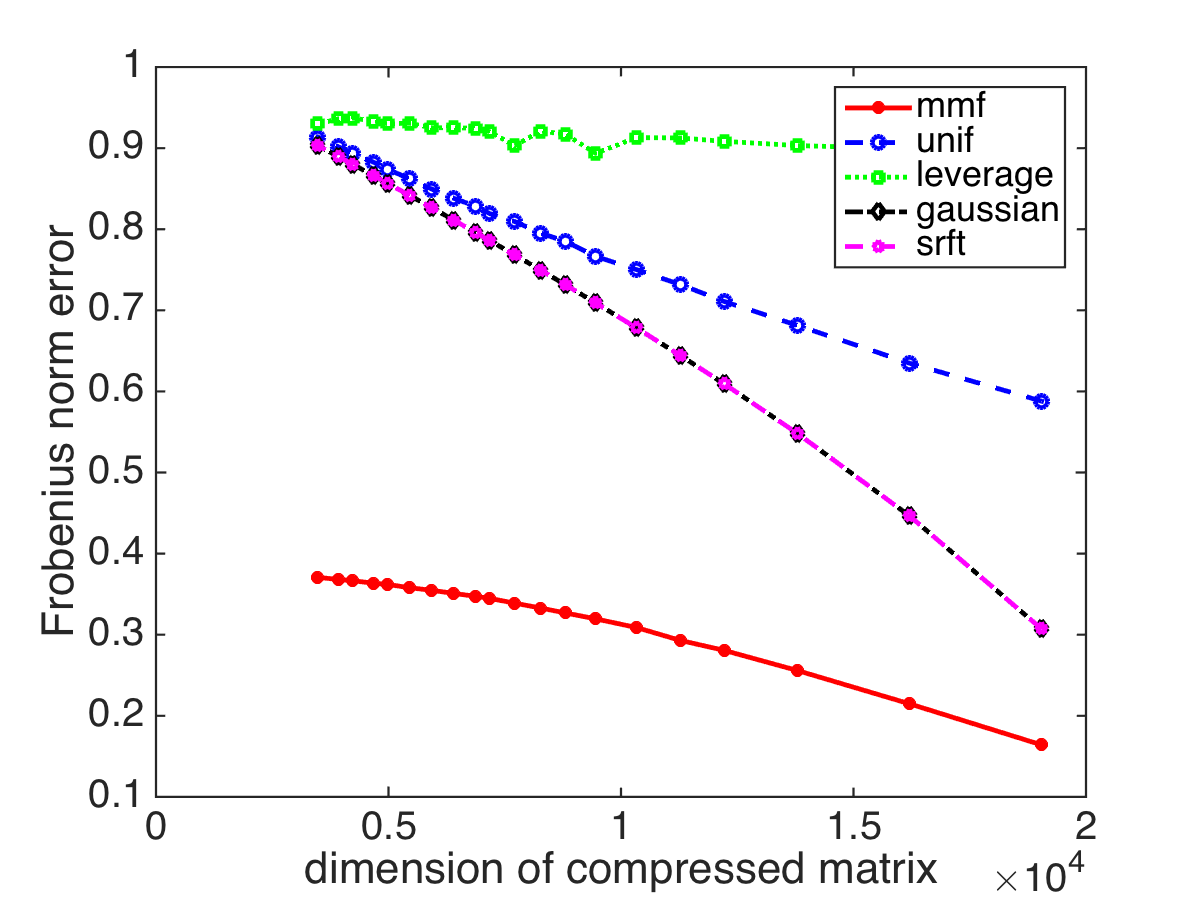}}
\end{minipage}
\begin{minipage}{.245\textwidth}
\centerline{\includegraphics[width=1.12\textwidth]{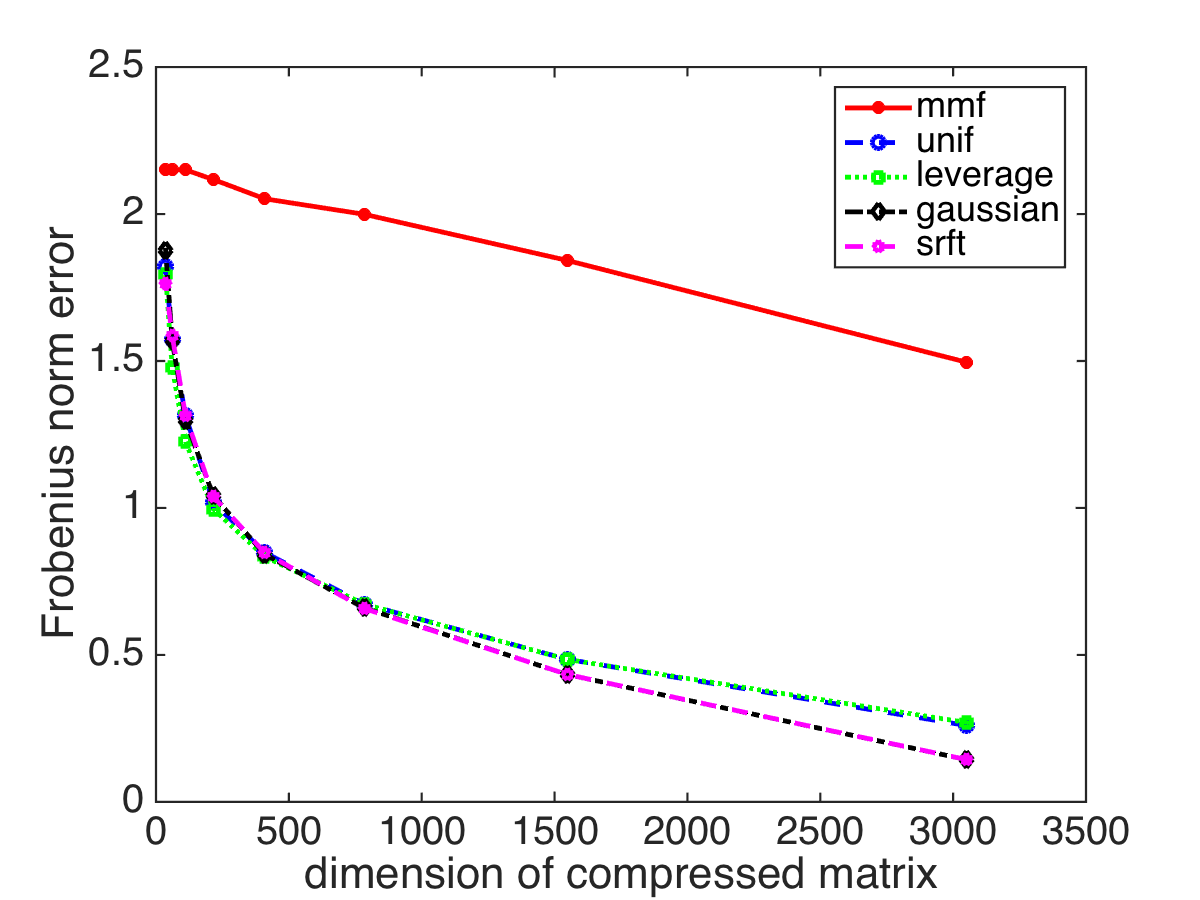}}
\end{minipage}
\vspace{-10pt}
\caption{\label{fig: compress}  
The Frobenius norm error \m{\nmN{A-\tilde A}_{\text{Frob}}} of compressing matrices with pMMF 
vs.\;other sketching methods, as a function of the dimension of the compressed core. In each figure, 
the error is normalized by \m{\nmN{A-A_k}_{\Frob}}, where \m{A_k} is the best rank \m{k} approximation 
to \m{A}. The four datasets are ``HEPph'', ``AstroPh'', ``CondMat'' and ``Gisette'', 
with \m{k\<=100} in the first three and \m{k\<=12} in ``Gisette''.}
\end{figure}

The main distinction between pMMF and other matrix sketching methods is that while the latter, 
implicitly or explicitly, make the assumption that \m{A} is low rank, or close to low rank, 
MMF makes a different structural assumption about \m{A}, namely that it has hidden 
hierarchical or multiresolution structure. 
Figure \ref{fig: compress} shows the results of experiments comparing the perfomance of MMF to 
other matrix sketching algorithms on some standard datasets. 
As in panes 1--3, on most datasets that we tried, pMMF significantly outperforms the 
other sketching methods in both Frobenius norm error and spectral norm error (plots of the latter 
can be found in the Supplement). The advantage of pMMF seems to be particularly great on network 
graphs, perhaps not surprisingly, since it has long been conjectured that networks have multiresolution 
structure \cite{Ravasz2003,Coifman06,Savas}. However, we find that pMMF often outperforms other methods on kernel matrices in general. 
On the other hand, on a small fraction of datasets, typically those which explicitly have low rank 
or are very close to being low rank (e.g., the fourth pane of Figure \ref{fig: compress}),   
pMMF performs much worse than expected. 
In such cases, a combination of the low rank and multiresolution approaches might be most 
advantageous, which is the subject of ongoing work. 

It is important to emphasize that pMMF is very scalable. Many other Nystr\"om methods are implemented in 
MATLAB, which limits the size of datasets on which they can be feasibly ran on. 
Moreover, leverage score methods require estimating the singular vectors of \m{A}, which, unless 
\m{A} is very low rank, can be a computational bottleneck. 
Several of the Nystr\"om experiments took 30 minutes or more to run on 8 cores, whereas 
our custom C++ pMMF implementation compressed the matrix in at most one or two minutes 
(see more timing results in the Supplement). 
Hence, pMMF addresses a different regime than many other Nystr\"om papers: 
whereas the latter often focus on compressing \m{\sim\!10^3} dimensional matrices to just 
10--100 dimensions, we are more interested in compressing \m{\sim\!10^4}--\m{10^5} dimensional 
matrices to \m{\sim\!10^3} dimensions.

\section{pMMF Preconditioning}\label{sec: preconditioning}

Given a large matrix \m{A\tin\RR^{n\times n}}, and a vector \m{b\tin\RR^n}, solving 
the linear system \m{Ax\<=b} is one of the most fundamental problems in computational mathematics. 
In the learning context, solving linear systems is the computational bottleneck in a range 
of algorithms. 
When \m{A} is sparse, \m{Ax\<=b} is typically solved with iterative methods, such as 
conjugate gradients. 
However, it is well known that the number of iterations needed for such methods to converge  
scales badly with \m{\kappa}, where \m{\kappa} is the ratio of the largest and smallest eigenvalues of \m{A}, 
called the \df{condition number}. 

The idea of \df{preconditioning} is to solve instead of \m{Ax\<=b} the related system \m{(M^{-1}\nts A)\ts x\<=M^{-1}b}, 
where \m{M^{-1}} is an easy to compute rough approximation to \m{A^{-1}}. 
A good preconditioner will ensure that 
\m{M^{-1}\nts A} is fast to multiply with the current vector iterate, 
while the condition number of \m{M^{-1}\nts A} is much better than that of the original matrix \m{A}. 
At the same time, it is important that \m{M^{-1}} be easily computable for massive matrices.  
For symmetric matrices, a variation on the above is to solve 
\m{(M^{-1/2} A\tts M^{-1/2})\ts y\<=M^{-1/2}\ts b}, and then set \m{x\<=M^{-1/2} y}, which retains symmetry. 
pMMF is a natural candidate preconditioner for symmetric matrices since (a) the pMMF factorization 
can be computed very fast (b) as evidenced by the previous section, \m{\tilde A} is a good 
approximation to \m{A}, (c) \m{\tilde A^{-\xi}} (with \m{\xi\tin\cbrN{1,1/2}}) can be computed from \m{\tilde A} 
in just \m{O(n+\delta_L^3)} time. However, unlike some other preconditioners, \m{\tilde A^{-\xi}} is 
generally not sparse. Therefore, in MMF preconditioning one never expands \m{\tilde A^{-\xi}} 
into a full matrix, but rather \m{\tilde A^{-\xi}} is applied to the vectors involved in the 
iterative method of choice as a sequence of rotations, as described in Section \ref{sec: sparsity}. 

\begin{figure}
\begin{minipage}{.24\textwidth}
\centerline{
\includegraphics[width=1.13\textwidth]{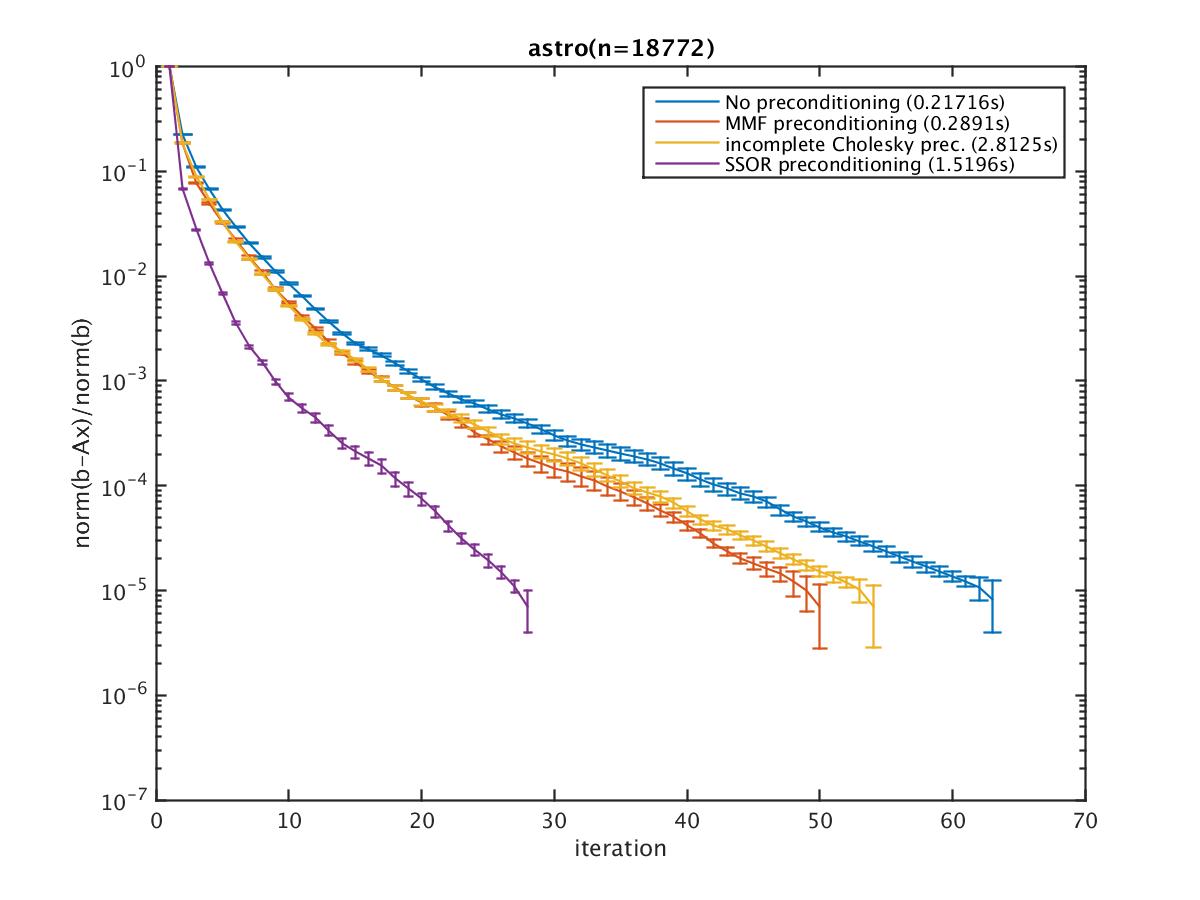}
}\end{minipage}
\begin{minipage}{.24\textwidth}
\centerline{\includegraphics[width=1.13\textwidth]{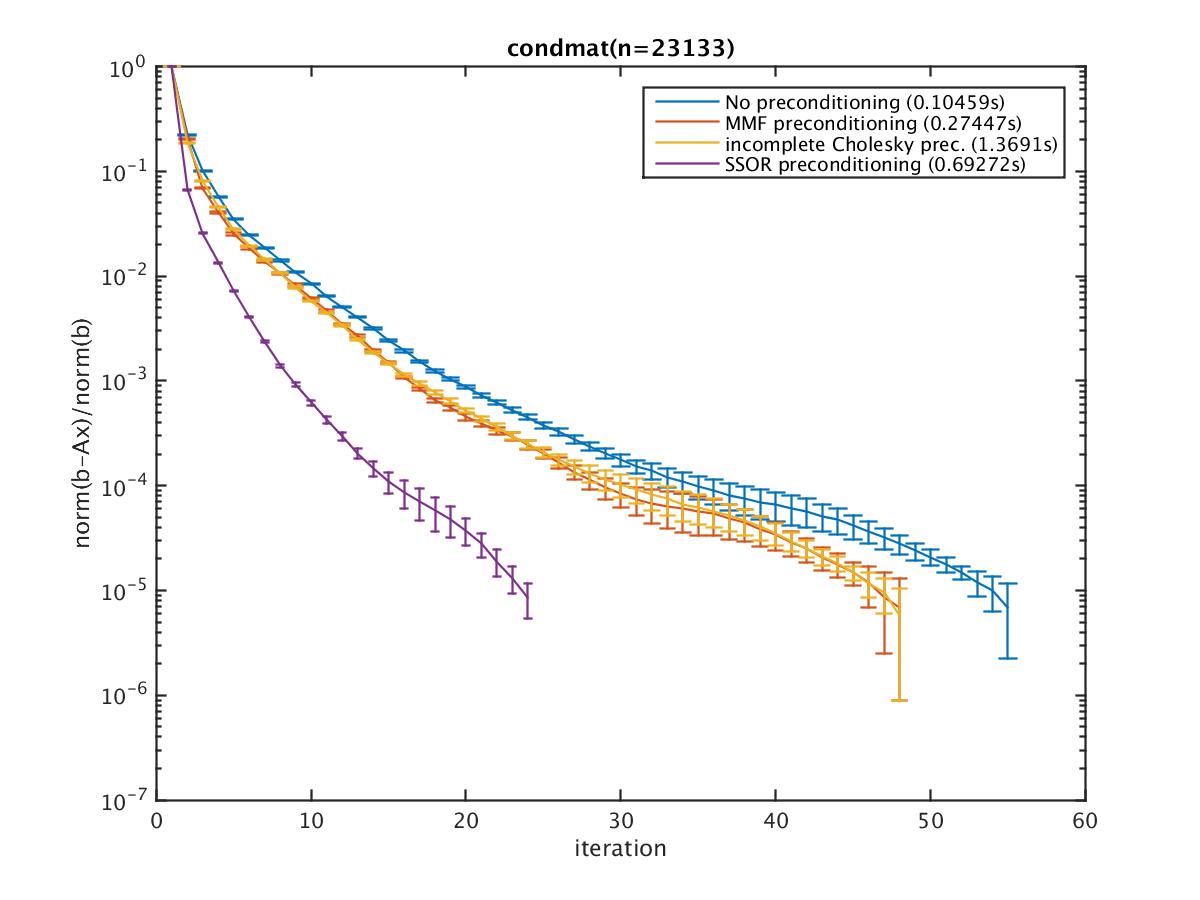}}
\end{minipage}
\begin{minipage}{.24\textwidth}
\centerline{\includegraphics[width=1.13\textwidth]{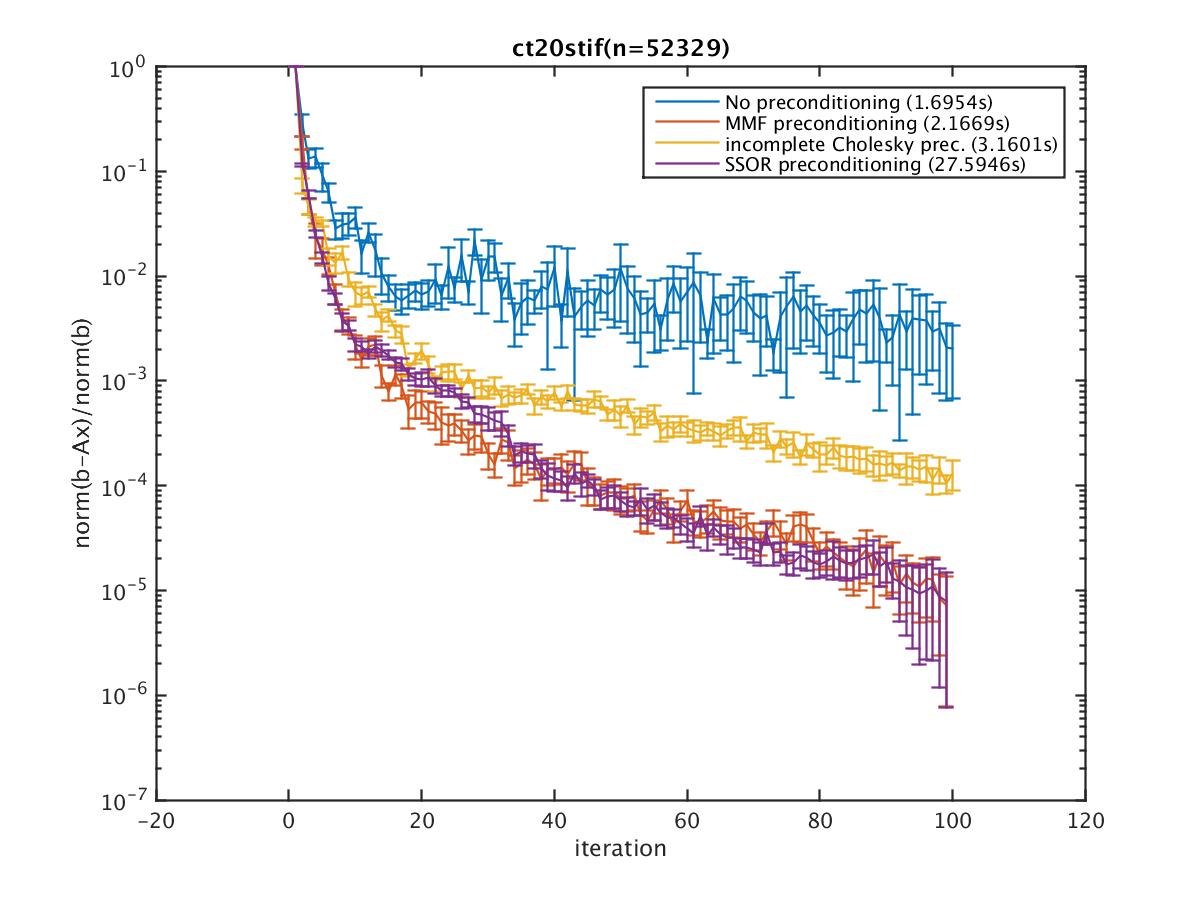}}
\end{minipage}
\begin{minipage}{.24\textwidth}
\centerline{\includegraphics[width=1.13\textwidth]{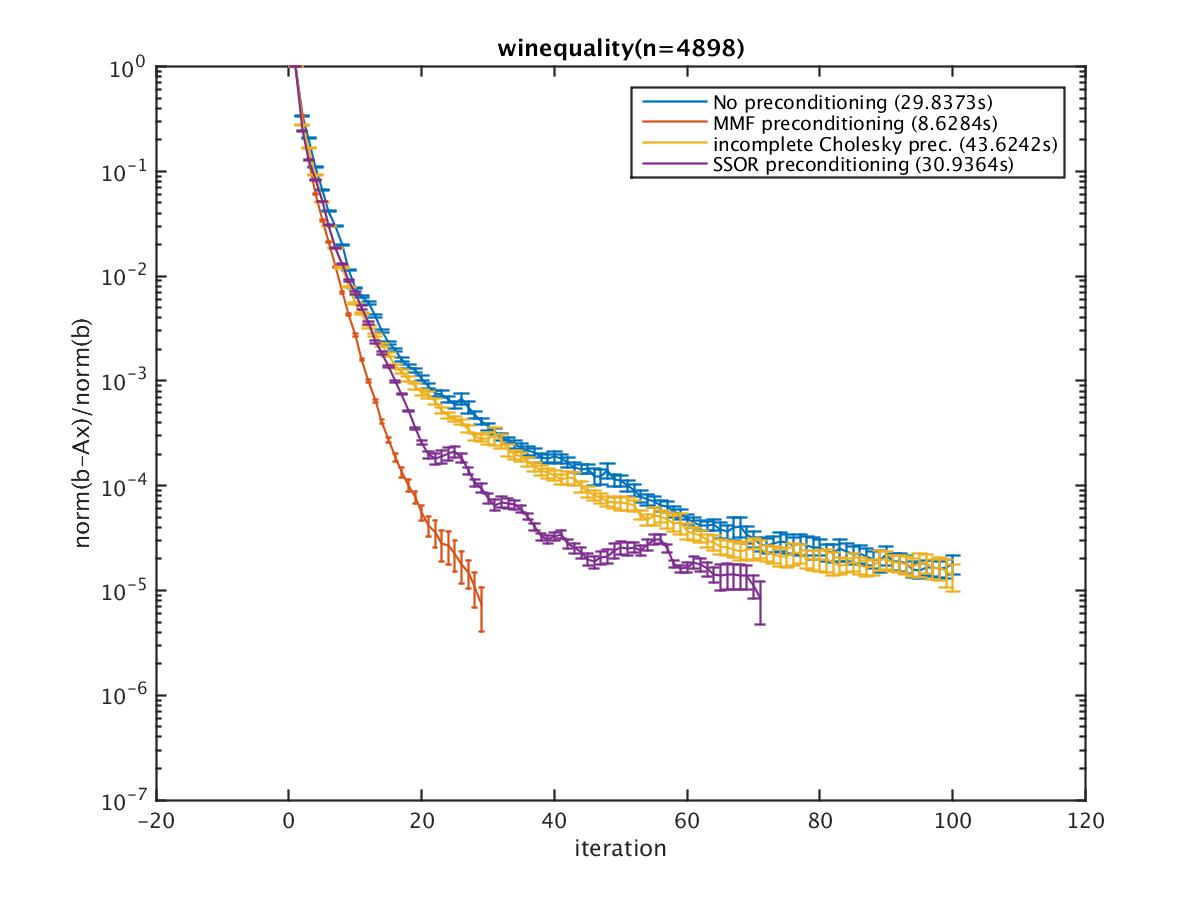}}
\end{minipage}
\vspace{-7pt}
\caption{\label{fig: linsolve} 
Residual as a function of conjugate gradient iterations when solving \m{Ax\<=b}, where \m{b} is a dense 
random vector. The indicated times are the wall clock time to convergence on an 8 core \m{2.6} GHz machine. 
The per-iteration time of pMMF preconditioning is usually 2--10 times faster than of other methods. 
pMMF is indicated in red. 
}
\end{figure}

A large number of different preconditioners have been proposed in  the literature, and  
even for a given type of problem there is often no single best choice, rather the choice reduces to 
experimentation. In our experiments our goal was to show that pMMF preconditioning 
can improve the convergence of linear solvers in learning problems, and that it is competitive 
with other preconditioners. Figure \ref{fig: linsolve} compares the performance of pMMF as a preconditioner 
to other standard preconditioners, such as incomplete Cholesky and SSOR. 
Several more preconditioning results are presented in the Supplement. 
In summary, pMMF appears competitive with other preconditioners on network and kernel matrices, and 
sometimes outperforms other methods. 
In each of our experiments, the time required to compute the pMMF preconditioner was less than a minute, 
which is amortized over the number of linear solves. 
We are still experimenting with 
how much pMMF as a preconditioner can be improved by fine tuning its parameters, and how well it will 
perform coupled with other solvers. 

\section{Conclusions}

The most common structural assumption about large matrices arising in learning problems 
is that they are low rank.
This paper explores the alternative approach of assuming that they have multiresolution structure. 
Our results suggest that not only is the multiresolution model often more faithful to the actual 
structure of data (e.g., as evidenced by much lower approximation error in compression experiments), 
but it also lends itself to devising efficient parallel algorithms, 
which is critical to  dealing with large scale problems. 
Our approach bears some similarities to multigrid methods \cite{Brandt73} 
and structured matrix decompositions \cite{Hackbusch1999,Borm2007,Chandrasekaran2005}, 
which are extremely popular in applied mathematics, primarily  
in the context of solving systems of partial differential equations. 
A crucial difference, however, is that whereas in these algorithms the multiresolution structure 
is suggested by the geometry of the domain, in learning problems the structure itself 
has to be learnt ``on the fly''. 
Empirically, the pMMF algorithm described in this paper scales linearly in the size of the data. 
Further work will explore folding entire learning and optimization algorithms into 
the multiresolution framework, while retaining the same scaling behavior. 
\newcommand{\nclust}{m^{\phantom{2}}} 
\newcommand{\nclustsq}{m^2} 
\begin{table}[h!]
\centerline{
\renewcommand{\arraystretch}{1.1}
\begin{tabular}{lccccccc}
\hline
& \multicolumn{2}{c}{serial MMF} & \multicolumn{2}{c}{pMMF operations}& \multicolumn{2}{c}{pMMF time}&\\
&dense&sparse&dense&sparse&dense&sparse&\m{N_{\text{proc}}}\\
\hline
Computing Grams&
\m{O(n^3)}&\m{O(\gamma n^3)}&
\m{O(pcn^2)}&\m{O(\gamma pcn^2)}&
\m{O(pc^3)}&\m{O(\gamma pc^3)}&
\m{\nclustsq}\\
Finding Rotations& 
\m{O(n^3)}&\m{O(n^3)}&
\m{O(cn)}& \m{O(cn)}&
\m{O(c^2)}&\m{O(c^2)}&
\m{\nclust}\\
Updating Grams& 
\m{O(n^3)}&\m{O(\gamma^2n^3)}&
\m{O(c^2n)}& \m{O(\gamma^2 c^2n)}&
\m{O(c^3)}&\m{O(\gamma^2c^3)}&
\m{\nclust}\\
Applying rotations &
\m{O(kn^2)}&\m{O(\gamma k n^2)}&
\m{O(kn^2)}&\m{O(\gamma k n^2)}&
\m{O(kc^2)}&\m{O(\gamma k c^2)}&
\m{\nclustsq}\\
Clustering&&&
\m{O(pmn^2)}&\m{O(\gamma pmn^2)}&
\m{O(pcn)}&\m{O(\gamma pcn)}&
\m{\nclustsq}\\
Reblocking&&&
\m{O(pn^2)}&\m{O(\gamma pn^2)}&
\m{O(pcn)}&\m{O(\gamma pcn)}&
\m{\nclust}\\
\hline
Factorization total&
\m{O(n^3)}&\m{O(n^3)}&
\m{O(pcn^2)}&\m{O(\gamma pcn^2)}&
\m{O(pc^3)}&\m{O(\gamma p c^3)}&
\m{\nclustsq}\\
\hline
\ignore{
Computing \m{\tilde A^{-1}}& 
\multicolumn{2}{c}{\m{O(g^3\nts\<+n)}}&
\multicolumn{2}{c}{\m{O(g^3\nts\<+n)}}&
\multicolumn{2}{c}{\m{O(g^3\nts\<+n)}}&
\m{1}\\
Computing \m{\tilde A v}&
\multicolumn{2}{c}{\m{O(kn)}}&
\multicolumn{2}{c}{\m{O(kn)}}&
\multicolumn{2}{c}{\m{O(kn)}}&
\m{1}\\
\hline}
\end{tabular}
}\vspace{-5pt}
\caption{\label{tbl: complexity} The rough complexity of different subtasks in 
pMMF vs.\;in the original serial MMF algorithm of \cite{MMFicml2014}.
Here \m{n} is the dimensionality of the original matrix, \m{A}, \m{k} 
is the order of the rotations, and \m{\gamma} is the fraction of non-zero entries in 
\m{A}, when \m{A} is sparse. We neglect that during the course of the computation 
\m{\gamma} tends to increase, because concomitantly \m{A_\ell} shrinks, and 
computation time is usually dominated by the first few stages. 
We also assume that entries of sparse matrices can be accessed in constant time. 
In pMMF, \m{p} is the number of stages, \m{m} is the number of 
clusters in each stage, and \m{c} is the typical cluster size (thus, \m{c\<=\theta(n/m)}). 
The ``pMMF time'' columns give the time complexity of the algorithm assuming 
an architecture that affords \m{N_{\text{proc}}}--fold parallelism. 
\m{g\<=\delta_L} is the size of the dense core in \m{H}.      
It is assumed that \m{k\<\leq p\<\leq c\<\leq n}, but \m{n\<=o(c^2)}.  
} 
\end{table}

\clearpage
{\small
\bibliographystyle{unsrt}
\bibliography{../bibliography/matrices,../bibliography/ml,../bibliography/graphs,../bibliography/harmonic,../bibliography/probability,../bibliography/numerical,pMMF}

\begin{thebibliography}{10}

\bibitem{ho2013more}
Qirong Ho, James Cipar, Henggang Cui, Seunghak Lee, Jin~Kyu Kim, Phillip~B
  Gibbons, Garth~A Gibson, Greg Ganger, and Eric~P Xing.
\newblock More effective distributed ml via a stale synchronous parallel
  parameter server.
\newblock In {\em Advances in neural information processing systems}, pages
  1223--1231, 2013.

\bibitem{parameter}
Mu~Li, Li~Zhou, Zichao Yang, Aaron Li, Fei Xia, David~G Andersen, and AJ~Smola.
\newblock Parameter server for distributed machine learning.
\newblock In {\em Big Learning NIPS Workshop}, 2013.

\bibitem{lee2014model}
Seunghak Lee, Jin~Kyu Kim, Xun Zheng, Qirong Ho, Garth~A Gibson, and Eric~P
  Xing.
\newblock On model parallelization and scheduling strategies for distributed
  machine learning.
\newblock In {\em Advances in Neural Information Processing Systems}, pages
  2834--2842, 2014.

\bibitem{MMFicml2014}
Risi Kondor, Nedelina Teneva, and Vikas Garg.
\newblock {M}ultiresolution {M}atrix {F}actorization.
\newblock In {\em Proceedings of the 31st International Conference on Machine
  Learning (ICML-14)}, pages 1620--1628, 2014.

\bibitem{Coifman2004}
Ronald~R. Coifman, Ronald~R Coifman, and Mauro Maggioni.
\newblock {Multiresolution analysis associated to diffusion semigroups:
  construction and fast algorithms}.
\newblock 2004.

\bibitem{Lee2008}
Ann~B Lee, Boaz Nadler, and Larry Wasserman.
\newblock {Treelets {---} An adaptive multi-scale basis for sparse unordered
  data}.
\newblock {\em Annals of Applied Statistics}, 2(2):435--471, 2008.

\bibitem{buluc2012}
Aydin Bulu{\c{c}} and John~R Gilbert.
\newblock Parallel sparse matrix-matrix multiplication and indexing:
  implementation and experiments.
\newblock {\em {SIAM} J Sci Comput}, 34(4), 2012.

\bibitem{Halko2011}
Nathan Halko, Per-Gunnar Martinsson, and Joel~A Tropp.
\newblock Finding structure with randomness: Probabilistic algorithms for
  constructing approximate matrix decompositions.
\newblock {\em SIAM review}, 53(2):217--288, 2011.

\bibitem{Ailon2009}
Nir Ailon and Bernard Chazelle.
\newblock {The Fast Johnson--Lindenstrauss Transform and Approximate Nearest
  Neighbors}.
\newblock {\em SIAM Journal on Computing}, 39(1):302--322, January 2009.

\bibitem{fastfood2013}
Quoc Le, Tam{\'a}s Sarl{\'o}s, and Alexander Smola.
\newblock Fastfood: -- computing {H}ilbert space expansions in loglinear time.
\newblock {\em JMLR}, 2013.

\bibitem{Williams2001}
C.~Williams and M.~Seeger.
\newblock Using the {N}ystr\"{o}m method to speed up kernel machines.
\newblock In {\em Advances in Neural Information Processing Systems (NIPS)},
  2001.

\bibitem{Fowlkes2004}
Charless Fowlkes, Serge Belongie, Fan Chung, and Jitendra Malik.
\newblock {Spectral grouping using the Nystr\"{o}m method.}
\newblock {\em IEEE transactions on pattern analysis and machine intelligence},
  26(2):214--25, February 2004.

\bibitem{Drineas2006all}
P.~Drineas, R.~Kannan, and M.~W. Mahoney.
\newblock {Fast monte carlo algorithms for matrices {I}--{III}}.
\newblock {\em SIAM Journal on Computing}, 36(1):158--183, 2006.

\bibitem{Mahoney2009CUR}
Michael~W Mahoney and Petros Drineas.
\newblock {CUR} matrix decompositions for improved data analysis.
\newblock {\em Proceedings of the National Academy of Sciences},
  106(3):697--702, 2009.

\bibitem{MahoneyRandomized}
M.~W. Mahoney.
\newblock {Randomized algorithms for matrices and data}.
\newblock {\em Foundations and Trends in Machine Learning}, 3, 2011.

\bibitem{Gittens2013}
Alex Gittens and Michael~W Mahoney.
\newblock {Revisiting the Nystr{\"o}m method for improved large-scale machine
  learning}.
\newblock volume~28, 2013.

\bibitem{Kumar2009}
Sanjiv Kumar, Mehryar Mohri, and Ameet Talwalkar.
\newblock {Ensemble Nystr\"{o}m Method}.
\newblock 2009.

\bibitem{Zhang2010}
Kai Zhang and James~T Kwok.
\newblock {Clustered Nystr\"{o}m method for large scale manifold learning and
  dimension reduction.}
\newblock {\em IEEE transactions on neural networks / a publication of the IEEE
  Neural Networks Council}, 21(10):1576--87, October 2010.

\bibitem{Deshpande2006}
Amit Deshpande, Luis Rademacher, Santosh Vempala, and Grant Wang.
\newblock Matrix approximation and projective clustering via volume sampling.
\newblock In {\em Proceedings of the seventeenth annual ACM-SIAM symposium on
  Discrete algorithm}, pages 1117--1126. ACM, 2006.

\bibitem{Kumar2012}
Sanjiv Kumar, Mehryar Mohri, and Ameet Talwalkar.
\newblock {Sampling Methods for the Nystr\"om method}.
\newblock {\em Journal of Machine Learning Research}, 13:981--1006, 2012.

\bibitem{Wang2013}
Shusen Wang and Zhihua Zhang.
\newblock {Improving CUR Matrix Decomposition and the Nystr\"om Approximation
  via Adaptive Sampling}.
\newblock 14:2729--2769, 2013.

\bibitem{Ravasz2003}
Erzs\'{e}bet Ravasz and Albert-L\'{a}szl\'{o} Barab\'{a}si.
\newblock {Hierarchical organization in complex networks}.
\newblock {\em Physical Review E}, 67(2):026112, February 2003.

\bibitem{Coifman06}
R~R Coifman and M~Maggioni.
\newblock {Diffusion wavelets}.
\newblock {\em Applied and Computational Harmonic Analysis}, 2006.

\bibitem{Savas}
Berkant Savas, Inderjit~S Dhillon, et~al.
\newblock Clustered low rank approximation of graphs in information science
  applications.
\newblock In {\em SDM}, pages 164--175. SIAM, 2011.

\bibitem{gpl}
The {GNU} {P}ublic {L}icense, {V}ersion 3, \url{http://www.gnu.org/licenses/}.

\bibitem{Brandt73}
Achi Brandt.
\newblock Multi-level adaptive technique (mlat) for fast numerical solution to
  boundary value problems.
\newblock In Henri Cabannes and Roger Temam, editors, {\em Proceedings of the
  Third International Conference on Numerical Methods in Fluid Mechanics},
  volume~18 of {\em Lecture Notes in Physics}, pages 82--89. Springer Berlin
  Heidelberg, 1973.

\bibitem{Hackbusch1999}
Wolfgang Hackbusch.
\newblock {A Sparse Matrix Arithmetic based on H -Matrices . Part I :
  Introduction to H -Matrices â}.
\newblock 62:1--12, 1999.

\bibitem{Borm2007}
Steffen Borm.
\newblock {Construction of data-sparse $H^2$-matrices by hierarchical
  compression}.
\newblock pages 1--33, 2007.

\bibitem{Chandrasekaran2005}
S.~Chandrasekaran, M.~Gu, and W.~Lyons.
\newblock A fast adaptive solver for hierarchically semiseparable
  representations.
\newblock {\em Calcolo}, 42(3-4):171{\textendash}185, 2005.

\end{thebibliography}
}

\end{document}